\documentclass[11pt]{amsart}

\usepackage{amsfonts,amsmath,amscd,latexsym,amssymb}

\input xy

\xyoption{all}

\usepackage{mathrsfs}

\theoremstyle{plain}

\newtheorem{teo}{Theorem}[section]
\newtheorem{propo}[teo]{Proposition}

\newtheorem{lema}[teo]{Lemma}

\theoremstyle{definition}

\theoremstyle{remark}

\newcommand{\Te}{\Theta}

\newcommand{\al}{\alpha}

\newcommand{\I}{\mathcal{I}}

\renewcommand{\O}{\mathcal{O}}

\newcommand{\ra}{\rightarrow}

\newcommand{\g}{\tilde{g}}
\newcommand{\upC}{\widetilde{C}}
\newcommand{\F}{\mathcal{F}}
\newcommand{\V}{\widetilde{V}}

\newcommand{\Z}{\mathbb{Z}}

\DeclareMathOperator{\codim}{codim}
\DeclareMathOperator{\Pic}{Pic}
\DeclareMathOperator{\Nm}{Nm}

\newcommand{\set}[1]{\left\{#1\right\}}
\newcommand{\abs}[1]{\left\vert#1\right\vert}
\newcommand{\apl}[5]{\left.\begin{array}{rrcl} #1:& \!\!\!\!#2&\longrightarrow& #3\\
&\!\!\!\! #4&\longmapsto& #5\end{array}\right.}
\newcommand{\apli}[4]{\left.\begin{array}{rcl} #1&\longrightarrow& #2\\
#3&\longmapsto& #4\end{array}\right.}
\newenvironment{sis}{\left\{\begin{aligned}}{\end{aligned}\right.}

\theoremstyle{plain}
\newtheorem{teoX}{Theorem}
\newtheorem{corX}[teoX]{Corollary}

\newtheorem*{deff*}{Definition}
\newtheorem*{teo*}{Theorem}
\newtheorem*{conj*}{Conjecture}
\newtheorem*{que*}{Question}

\title[Cohomological support loci for Abel-Prym
curves]{Cohomological support loci for Abel-Prym curves}

\author[S. Casalaina-Martin]{Sebastian Casalaina-Martin}
\thanks{The first author was partially supported by NSF MSPRF grant DMS-0503228}
\address{
Harvard University,
Department of Mathematics,
One Oxford Street,
Cambridge, MA 02138
USA }
\email{casa@math.harvard.edu}

\author[M.~Lahoz]{Mart\'{i} Lahoz}
\thanks{The second author has been supported by  Ministerio de Educaci\'on y
Ciencia, beca de Formaci\'on del Profesorado Universitario,
MTM2006-14234-C02 and 2005SGR-00557.}
\address{ Departament de Matem\`atica  Aplicada I, Universitat Polit\`ecnica de
Cata-lunya,
08028 Barcelona (Spain).}
\email{{\tt marti.lahoz@upc.edu}}

\author[F.~Viviani]{Filippo Viviani}
\address{Institut f\"ur Mathematik, Humboldt Universit\"at zu Berlin, 10099 Berlin (Germany).}
\email{{\tt viviani@math.hu-berlin.de}}

\date{\today}

\thanks{2000\,\emph{Mathematics Subject Classification:} 14H40, 14K12, 14F17}

\begin{document}

\begin{abstract}
For  an Abel-Prym curve contained in a Prym variety, we determine the cohomological
support loci of its twisted ideal sheaves and the dimension of its theta-dual.
\end{abstract}

\maketitle

\section*{Introduction}

The purpose of this paper is to study the theta-dual and the
cohomological support loci for the twisted ideal sheaves of an
Abel-Prym curve contained in the Prym variety
associated to an \'etale double cover of smooth projective non-hyperelliptic
curves.


Recall that given a coherent sheaf $\F$ on a smooth
projective variety $X$, the $i$-th cohomological support locus of
$\F$ is
$$
V^i(\F):=\{\alpha\in \Pic^0(X)\: |\: h^i(X, \F\otimes \alpha)>0\}\subset \Pic^0(X).
$$
These loci have been studied in a number of contexts, and were considered for example by
Green-Lazarsfeld (see \cite{GL87,GL91}) in order to prove a
generic vanishing theorem for the canonical sheaf of irregular
complex varieties. More precisely, they proved that if the
Albanese morphism $a: X\to {\rm Alb}(X)\cong\widehat{\Pic^0(X)}$
has generic fiber of dimension $k$, then
\begin{equation*}
\codim_{\Pic^0(X)} V^i(\omega_X)\geq i-k.
\end{equation*}

In a sequence of articles (\cite{P,PP06,PP07} cf. \cite{PP05})), G.~Pareschi and M.~Popa have studied similar questions, and
have introduced the notion of a $GV_k$ sheaf (see also \cite{Hac}):
\begin{deff*}
A coherent sheaf $\F$ is said to be $GV_k$ for some $k\in \Z$
(which stands for generic vanishing of order $k$) if
\begin{equation*}
\codim_{\Pic^0(X)} V^i(\F)\geq i-k \: \text{ for all } i>0.
\tag{*}
\end{equation*}
\end{deff*}
Observe that we have the natural inclusions between the
$GV_k$-sheaves:
$$(W)IT_0=GV_{-d}\subset \cdots \subset GV_{-1}=M\subset GV_0=GV\subset \cdots \subset
GV_d={\rm Coh}(\F)
\footnote{In the above
chains of inclusions we have also indicated some other names that are
used in the literature, namely: the sheaves satisfying $GV_0$ are also
called $GV$-sheaves (which stands for generic vanishing sheaves), the
$GV_{-1}$-sheaves  are called $M$-regular sheaves (which stands Mukai
regular sheaves) and the $GV_{-d}$-sheaves  are called $(W)IT_0$-sheaves
because they satisfy $V^i(\F)=\emptyset$ for all $i>0$, or in
other words they satisfies the (weak) index theorem with index $0$
in  Mukai's terminology.}
,$$ where the sequence becomes stationary
outside the interval $[-d, d]$ for $d=\dim(X)$.  With this terminology, the result of M.~Green and
R.~Lazarsfeld says that if the Albanese morphism has generic fiber of
dimension $k$, then $\omega_X$ is $GV_k$.


The above condition
$(*)$ can be expressed in terms of the Fourier-Mukai transform with respect to the
kernel $\mathcal{P}=(a, {\rm id})^*(\mathcal{L})\in \mathbb{D}(X\times
\hat{A})$, where $\mathcal{L}$ is the Poincar\'e line bundle on
$A\times \hat{A}$ (suitably normalized) and $A={\rm Alb}(X)$.  G.~Pareschi and M.~Popa use this to study  vanishing in a variety of contexts, and in particular  when $\F$ is an adjoint linear series or
an ideal sheaf suitably twisted; this has produced many
interesting applications (see \cite{PP1, PP2,PP3, PP06, PP07}).

When $X=A$ is an abelian variety, the case we will be considering here, $\Pic^0(A)$ is canonically
isomorphic to $\widehat A$, and so $V^i(\F)$ can be considered as a subspace
inside the dual abelian variety. Moreover, in the case of a  principally polarized abelian variety
$(A,\Te)$, which we will abbreviate by  ppav in the sequel, we can view the
cohomological support loci as subspaces of the original abelian
variety using the isomorphism $\varphi_\Te:A\stackrel\sim\ra \hat
A$. In this setting, Pareschi-Popa introduced the following definition of
theta-dual (see  \cite[Def. 4.2]{PP07}), which can be viewed as the cohomological support locus of a twisted ideal sheaf.

\begin{deff*}
Given a closed subscheme $X$ of a ppav $(A, \Theta)$, the dual
theta-dual $V(X)$ of $X$ is the closed subset\footnote{It is
possible to put a canonical schematic structure on $V(X)$ (see
\cite[Def. 4.2]{PP07}), which however will never play a role in
this paper.} defined by:
$$V( X):= V^0(\I_X(\Te)) = \{ \alpha\in \hat{A} ~|~ h^0 (A, \I_X( \Theta) \otimes \alpha) \neq 0\}
\subset \widehat A\cong \Pic^0(A).$$
\end{deff*}
\noindent  Observe that, via the principal polarization $\Theta$, the theta-dual $V(X)$ can be
canonically identified with the locus $\set{ a \in A ~|~ X \subset \Theta_a =t_a^*\Theta}$
of theta-translates containing $X$ (here $t_a:A\to A$ is the translation map).

More generally, we will be interested in the  cohomological support loci of the form $V^i(\I_X(n\Theta))$.
In particular, in this paper we consider the cohomological support loci associated to
the ideal sheaves of Abel-Prym curves.  To fix some notation, let $P$ be the Prym variety of
dimension $g-1$ associated to
the  \'etale double cover $\upC\to C$ of irreducible
smooth projective non-hyperelliptic curves of genus $\g=2g-1$ and $g\geq 3$,
respectively.  Let $\Xi$ be the canonical principal polarization.
Since $\upC$ is not hyperelliptic,
there is an embedding $\upC \hookrightarrow P$, unique up to
translation. We denote by $\I_{\upC}$ the ideal sheaf of $\upC$
inside $P$ (see section \ref{sec:notations} for the definition and
the standard notations).  Our main result is (combining Theorems \ref{dim},
\ref{1-twist}, \ref{2-twist}):

\begin{teoX}\label{main-theo}
Let $\upC$ be a non-hyperelliptic Abel-Prym curve embedded in its Prym variety $(P,\Xi)$. Then we have:
\begin{enumerate}
\item The theta-dual of $\upC$ has dimension
$$\dim (V(\upC))=\dim P-3=g-4.$$
 \item
The cohomological support loci of
$\I_{\upC}(\Xi)$ can be (non-canonically) identified with
$$\begin{sis}
&V^0(\I_{\upC} (\Xi))=V^1(\I_{\upC} (\Xi))=V(\upC), \\
&V^2(\I_{\upC} (\Xi))=P,\\
&V^{\geq 3}(\I_{\upC} (\Xi))=\emptyset.
\end{sis}
$$
\item The cohomological support loci of
$\I_{\upC}(2\Xi)$ can be (non-canonically) identified with
$$
\begin{sis}
&V^0(\I_{\upC} (2\Xi))=
\begin{cases}
P & \text{ if } g\geq 4,\\
\{0\}& \text{ if } g=3,
\end{cases}\\
&V^1(\I_{\upC} (2\Xi))=V^2(\I_{\upC} (2\Xi))=
\{0\},\\
&V^{\geq 3}(\I_{\upC} (2\Xi))=\emptyset.
\end{sis}$$
\end{enumerate}
\end{teoX}
The above theorem follows by combining Theorems \ref{dim}, \ref{1-twist},
\ref{2-twist} together with the identifications that we make in formulas
(\ref{auxiliary}) and (\ref{canonical-support}) below.
As an immediate Corollary we obtain:

\begin{corX}\label{main-cor}
With the notation of the above Theorem, the following hold
\begin{enumerate}
 \item[(i)] $\I_{\upC}(\Xi)$ is not $GV$.  More precisely, it is $GV_2$ but not $GV_1$;
 \item[(ii)] $\I_{\upC}(2 \Xi)$ is $GV$, but not $IT_0$.  More precisely, it is $GV_{-(g-2)}$ but not $GV_{-(g-1)}$;
 \item[(iii)] $\I_{\upC}(m\Xi)$ is $IT_0$ for every $m\geq 3$.
\end{enumerate}
\end{corX}

Statements (i) and (ii) follow by direct inspection from the theorem, while
part (iii) follows from part (ii) and the general fact
that if $\F$ is a $GV$-sheaf then $\F(\Theta)$ is $IT_0$ (see
\cite[Lemma 3.1]{PP07}).

These results should be compared to the case of an Abel-Jacobi curve.  Recall that for any non-rational curve $C$, the Abel-Jacobi map gives a non-canonical
embedding of $C$ in its canonically polarized
Jacobian $(J(C),\Te)$. In this case, the cohomological
support loci can be (non-canonically) identified with (combine \cite[Thm. 4.1, Prop. 4.4]{PP1}
and \cite[Lem. 3.3, Exa. 4.5]{PP07})
$$\begin{sis}
&V(C)=V^0(\I_C(\Te))=V^1(\I_C(\Te))=W_{g-2},\\
& V^0(\I_C(n\Te))=J(C) \text{ for } n\geq 2,\\
& V^{\geq 2}(\I_C(\Te))=V^{\geq 1}(\I_C(n\Te))=\emptyset \text{ for } n\geq 2,
\end{sis}$$
where $W_{g-2}=W_{g-2}^0$ is the Brill-Noether locus of line bundles of degree
$g-2$ with non-trivial global sections. From the above description, we get  that for an Abel-Jacobi curve
\begin{enumerate}
 \item[(i)] $\I_C(\Theta)$ is $GV$, but not $IT_0$.
 \item[(ii)] $\I_C(m\Theta)$ is $IT_0$ for every $m\geq 2$.
 \item[(iii)] ${\rm dim}(V(C))=g-2=\dim JC-2.$
\end{enumerate}

In \cite{PP07}, Pareschi-Popa have proved that the above condition
(i) characterizes Abel-Jacobi curves among the non-degenerate
curves inside a ppav. Moreover, they have further conjectured that the conditions (ii) and (iii) should also provide new
characterizations of Abel-Jacobi curves\footnote{We refer to
\cite{PP07} for analogous conjectures concerning the subvarieties
of a ppav of minimal cohomological class.}.

The results in this paper show that this conjecture is not violated by Abel-Prym curves,
which in a sense are the curves inside a ppav closest to the Abel-Jacobi
curves.

In addition, from the results on Abel-Prym curves above, and those cited for Abel-Jacobi curves,  it seems natural to ask for the relation between the following conditions on a curve $X$ on a ppav $(A,\Theta)$ of dimension $g$:

\begin{enumerate}

\item $\I_X(e\Theta)$ is $GV$, but not $IT_0$.
\item $\I_X((e+1)\Theta)$ is $IT_0$, but $\I_X(e\Theta)$ is not.
\item $\dim V(X)=g-e-1$.
\item $X$ is an Abel-Prym-Tyurin curve with Prym-Tyurin variety $(A, \Theta)\subset (JX, \Theta_X)$
of exponent $e$, that is $[X]=\frac{e[\Theta]^{g-1}}{(g-1)!}$.
\end{enumerate}

Prym-Tyurin varieties of exponent $1$ are precisely the Jacobians
(by the Matsusaka-Ran criterion), so that the Pareschi-Popa
conjecture states that the above conditions are all equivalent for $e=1$. In
the next case, it is known  (see \cite{W87}) that the closure
inside the moduli space $\mathcal A_g$ of ppav's of dimension $g$ of the
Prym-Tyurin varieties of exponent $2$ has a unique component of
maximal dimension (which is $3g$), namely the closure of the
classical Prym varieties\footnote{Among the Prym-Tyurin ppav of
dimension $g$ and exponent $2$, the classical Prym varieties can
also be characterized as those for which the curve $X$ is smooth
of maximal arithmetic genus, namely $2g+1$.}. Therefore, our
results in Theorem \ref{main-theo} show that ``most'' of the
Abel-Prym-Tyurin curves of exponent $2$ satisfies the conditions
(1), (2) and (3).

On the other hand, Andreas Hoering has pointed out to us that 
condition (3) is much weaker than condition (4):  any
curve $X$ on a subvariety $Y$ with $\dim(Y)>1$ and $\dim V(Y)=g-e-1$ will
have $\dim V(X)\ge g-e-1$. Since $Y$ will contain curves of arbitrarily high
degree with respect to $\Theta$, one can construct curves satisfying (3) but
not (4).  As a concrete example, consider a curve $X$ lying on a $W_d$
($1<d<e+1$) inside a Jacobian or on the Fano surface inside the intermediate
Jacobian of a cubic threefold.
Thus we propose an alternate version of (3), which may be more closely
related to the other conditions

\begin{enumerate}
\item[(3$'$)] $\dim V(X)=g-e-1$ and $X$ is not contained in a subvariety $Y$ with
$1<\dim(Y)<e+1$ and $\dim V(Y)=g-e-1$.
\end{enumerate}

Since an Abel-Prym curve of the intermediate Jacobian of a cubic threefold
lies on the Fano surface $F$, which has class $[\Theta]^3/3!$ and $\dim
V(F)=g-3$, (4) does not imply (3$'$), and so we suggest the following
modification of (4) as well:

\begin{enumerate}
\item[(4$'$)] $X$ is an Abel-Prym-Tyurin curve of exponent $e$ and $X$ is not contained in 
a subvariety $Y$ with $1<\dim(Y)=d<e+1$ and class $[Y]=\alpha \frac{[\Theta]^{g-d}}{(g-d)!}$ with
$\alpha<e$.
\end{enumerate}


The paper is organized as it follows. In section $2$, we review
the definition of the Prym variety $(P, \Xi)$ associated to an
\'etale double cover $\upC\to C$ in order to fix the notation
used throughout the paper. In section $3$, we
prove that the theta-dual of an Abel-Prym curve $\upC$ inside $P$
can be set-theoretically identified with the Brill-Noether locus
$V^2$ defined in \cite{W85}. General results about these
Brill-Noether (\cite{B,DP}) give the inequality $\dim V(\upC)\ge
\dim(P)-3$. Using ideas from \cite[sections 6, 7]{M}, we show that
the equality holds (Theorem \ref{dim}), which proves part $(1)$ of
the Main Theorem \ref{main-theo}.
 In section $4$ and $5$, we compute the cohomological support loci
for the twisted ideal sheaves $\I_{\upC}(\Xi)$ (Theorem
\ref{1-twist}) and $\I_{\upC}(2 \Xi)$ (Theorem \ref{2-twist}),
proving explicitly parts (2) and (3) of the Main Theorem
\ref{main-theo}.

\vspace{0,4cm}

\noindent\emph{Acknowledgments.} We are grateful to prof.
A.~Ragusa for organizing an excellent summer school ``Pragmatic
2007'' held at the University of Catania, where the three authors
began their joint work on this subject. We also thank prof. G.~Pareschi
and prof. M.~Popa who, during that summer school, suggested to us an interesting research problem from which this work originated,
and have since then followed the progress of this work, providing useful
suggestions.

\section{Notation and basic definitions}\label{sec:notations}

Throughout this paper, we work over an algebraically closed field
$k$ of characteristic different from $2$.
The basic results cited here are due to Mumford \cite{M71}.
Let $\pi: \upC
\rightarrow C$ be an \'etale double cover of irreducible
smooth projective curves of genus $\g$ and $g$, respectively. By
the Hurwitz formula, we get that $\g=2g-1$. We denote by
$\sigma$ the involution on $\upC$ associated to the above double
cover. Consider the norm map
$$\apli{\Nm : \Pic(\upC)}{\Pic(C)}{\O_{\upC}(\sum_j r_j p_j)}{\O_{C}(\sum_j r_j \pi(p_j)).}$$
The kernel of the norm map has two connected components
$$\ker \Nm= P\cup P'\subset \Pic^0(\upC),$$
where $P$ is the component containing the identity element and is, by definition,
the Prym variety associated to the \'etale double cover $\pi$. The
above components $P$ and $P'$ have the following explicit description
$$P=\set{\O_{\upC}(D-\sigma(D))\: | \: D\in {\rm Div}^{2N}(\upC), N\geq 0},$$
$$P'=\set{\O_{\upC}(D-\sigma(D))\: | \: D\in {\rm Div}^{2N+1}(\upC), N\geq 0}.$$
It is often useful to consider the inverse image of the canonical line bundle of $C$ via the norm map.  This also has
two connected components
$$\Nm^{-1}(\omega_C)=P^+\cup P^-\subset \Pic^{2g-2}(\upC)=\Pic^{\g-1}(\upC),$$
which have the following explicit description
$$P^+=\set{L\in \Nm^{-1}(\omega_C)\: | \: h^0(L)\equiv 0 \mod 2},$$
$$P^-=\set{L\in \Nm^{-1}(\omega_C)\: | \: h^0(L)\equiv 1 \mod 2}.$$
The above varieties $P'$, $P^+$ and $P^-$ are isomorphic to
the Prym variety $P$ and, in this work, we will pass frequently
from one to another.

\vspace{0,4cm}

There is a principal polarization $\Xi\in {\rm NS}(P)$ induced by the principal polarization $\Theta_{\upC}\in {NS}(J\upC)$.  In fact, $\Theta_{\upC}|_P=2\Xi$.
One of the primary motivations for considering $P^+$ is the existence of a canonically defined divisor $\Xi^+$ whose class in the Neron-Severi group of $P$
is $\Xi$:
$$
\Xi^{+}=\set{L\in P^+\subset \Pic^{\tilde g-1}(\upC)\mid h^0(L)>0 }\subset P^+.
$$

\noindent On the other hand, the canonical Abel-Prym map is defined as
$$
\apli{i: \upC}{P'}{p}{\sigma(p)-p.}
$$
If $\upC$ is hyperelliptic then the image of $\upC$
via the Abel-Prym map is a smooth hyperelliptic curve $D$ and
the Prym variety $P$ is isomorphic to the Jacobian $J(D)$ of $D$ (\cite[Cor.
12.5.7]{BL}).  On the other hand, if $C$ is hyperelliptic but $\upC$ is not, then
the Prym variety $P$ is the product of two hyperelliptic Jacobians (see \cite{M}).
Therefore, since we are mostly interested in the case of an irreducible
non-Jacobian ppav, we will assume throughout this paper that $C$ is not hyperelliptic
(and in particular $g\geq 3$).
Note that under this hypothesis, the
Abel-Prym map is an embedding (\cite[Cor. 12.5.6]{BL}).

\vspace{0,4 cm}

Since the Abel-Prym curve $\upC\subset P'$ and the polarization
$\Xi^+\subset P^+$ lie canonically in different spaces, the cohomological support
loci for the twisted ideal sheaf
$\I_{\upC}(n\Xi^+)$ is only defined up to a translation since we have to
choose a way to translate $\upC$ and $\Xi^+$ inside the Prym variety $P$.
For this reason, we introduce the following auxiliary (canonically
defined) loci
\begin{equation}\label{auxiliary}
\V^i(\I_{\upC}(n \Xi^+))=\{E\in P^-\: |\: h^i(P', \I_{\upC}(n \Xi^+_{E})>0\}\subset P^-,
\end{equation}
where $\I_{\upC}$ is the ideal sheaf of $\upC$ inside $P'$ and for
$E\in P^-\subset \Pic^{\g-1}(\upC)$, we denote by $\Xi^+_{E}\subset P'$ the
translate of the canonical theta divisor $\Xi^+$ by $E^{-1}$.
The relation between
$\V^i(\I_{\upC}(n \Xi^+))$ and $V^i(\I_{\upC}(n \Xi))$ is easy to work out.
In fact there is a choice of translate of $\upC\subset P$ and $\Xi\subset P$ so that under the isomorphism
\begin{equation*}
\apl{\psi_{E_0}}{P^-}{P}{E}{E\otimes E_0^{-1}}
\end{equation*}
induced by a  line bundle $E_0\in P^-$, we have
\begin{equation*}
V^i(\I_{\upC}(n\Xi))=\V^i(\I_{\upC}(n\Xi^+))^{\otimes n}\otimes E_0^{\otimes -n}:=\{(E\otimes E_0^{-1})^{\otimes n}\: |\: E\in \V^i(\I_{\upC}(n \Xi^+))\}.
\end{equation*}

For this reason, we will also identify the cohomological support loci with the following
canonical loci:
\begin{equation}\label{canonical-support}
V^i(\I_{\upC}(n\Xi))=\V^i(\I_{\upC}(n\Xi^+))^{\otimes n}\subset {\rm Nm}^{-1}
(\omega_C^{\otimes n})\subset \Pic^{n(\g-1)}(\upC).
\end{equation}


\vspace{0,4cm}

For later use, we end this section with the following Lemma, which describes
the restriction of the translates of the theta-divisor to the Abel-Prym curve.

\begin{lema}\label{restriction}
Given any $E\in P^-$, there is an isomorphism of line bundles
$$\O_{P'}(\Xi^+_{E})_{|\upC}\cong E.$$
Moreover, if $E\in P^- -  V(\upC)$ and $D$ is the unique divisor in $|E|$,
then we have an equality of divisors
$$\left(\Xi^+_{E}\right)_{\mid\upC}=\upC\cap \Xi^+_{E}= D.$$
\end{lema}
\begin{proof}
This is standard, we include
a proof for the convenience of the reader. Suppose first that
$E\in P^- - V(\upC)$, which, by Lemma \ref{dual-BN}, is equivalent
to the condition $h^0(\upC, E)=1$.
Write $\abs{E}=D=p_1+\ldots+p_{\g-1}$, where $p_i\in \upC$. Since $p_i$ is a fixed
point of the linear series $|E|$, we have that  $h^0(\upC, E\otimes \O_{\upC}(-p_i+\sigma p_i))=2$,
which implies that
$$D\subset \upC\cap \Xi^+_{E}=\left(\Xi^+_{E}\right)_{|\upC}.$$
Using that $\upC\cdot \Xi^+_{E}=\g-1$, we get the desired second equality.
Now consider the maps
$$ \upC \times P^- \stackrel{(a, {\rm id})}{\longrightarrow} P'\times P^-
\stackrel{\mu}{\longrightarrow} P^+,$$
where $a$ is the Abel-Prym map and $\mu$ is the multiplication map.
Let $\mathcal{P}$ be the Poincar\'e line bundle on $\upC\times P^-$, trivialized over the
section $\{p\}\times P^-$ for some $p\in \upC$.
Consider the line bundle on $\upC\times P^-$ given by $\mathcal{L}:=(a \times {\rm id})^*(\mu^*\O_{P^+}(\Xi^+))$. We can trivialize $\mathcal{L}$
along the given section $\{p\}\times P^-$ by tensoring with the pull back from $P^-$ of the divisor
$\Xi^+_{p-\sigma(p)}$.  It is easy to check that the fibers of $\mathcal{P}$ and $\mathcal{L}$
over $\upC\times \{E\}$ are given by
$$\begin{sis}
&\mathcal{P}_{\upC\times \{E\}}=E, \\
&\mathcal{L}_{\upC\times \{E\}}=\O_{P'}(\Xi^+_{E})_{|\upC}. \\
  \end{sis}$$
By what was proved above, if $E\not\in   P^- - V(\upC)$ then the two fibers agree.
By the Seesaw theorem (e.g. \cite[Lemma 11.3.4]{BL}), $\mathcal{P}\cong\mathcal{L}$
and we get the desired first equality.
\end{proof}

\section{The theta-dual of $\upC$}

In this section, we want to study the theta-dual of $\upC$ in
the Prym variety $P$; this can be identified canonically with the set
(see (\ref{canonical-support})):
$$V(\upC):=\set{E\in P^- \: | \: h^0(P', \I_{\upC}(\Xi^+_{E}))> 0}\subset P^-.$$

In fact the theta-dual $V(\upC)$ can be described in terms of the
following standard Brill-Noether loci (see \cite{W85}):
$$V^r:=\set{L\in {\rm Nm}^{-1}(\omega_{C})\: |\: h^0(L)\geq r+1,\: h^0(L)\equiv r+1\mod 2},$$
where $V^r\subset P^-$ (resp. $V^r\subset P^-$) if $r$ is even
(resp. odd). We view both the theta dual and the Brill-Noether
loci as sets, although they can be endowed with natural scheme
structures\footnote{We remark however that the $V^r$ have been
considered with different natural schematic structures (compare
\cite{W85} and \cite{DP}).}.

\begin{lema}\label{dual-BN}
We have the set-theoretic equality
$$V(\upC)=V^2.$$
\end{lema}
\begin{proof}
An element $E\in P^-$ belongs to $V(\upC)$ if and only if $\upC\subset
\Xi^+_{E}$, which, by the definition of $\upC\subset P'$,
is equivalent to $h^0(\upC, E\otimes \O_{\upC}(\sigma(p)-p))>0$ for every $p\in \upC$.
By  Mumford's parity trick (see \cite{M}), this happens if and only if  $h^0(\upC,E)\geq 3$, that is
$E\in V^2$.
\end{proof}

\begin{teo}\label{dim}
For any \'etale double cover $\upC\to C$ as above with $C$ non-hyperelliptic of genus $g$,
it holds that
$${\rm dim}(V^2)={\rm dim}(P)-3=g-4.$$
\end{teo}
\noindent For $g= 3$, the Theorem says that $V^2=\emptyset$.
We start with the following Lemma, which is similar to \cite[Lemma p. 345]{M}.

\begin{lema}\label{big-comp}
If $Z\subseteq V^2$ is an irreducible component, and
$\dim Z\ge g-3$, then for a general line bundle $L\in Z$, there is
a line bundle $M$ on $C$ with $h^0(M)\ge 2$, and an effective
divisor $F$ on $\upC$ such that $L\cong\pi^*M\otimes \mathscr
O_{\upC}(F)$.
\end{lema}
\begin{proof}
Let $Z$ and $L$ be as in the statement. Suppose that $h^0(L)=r+1$
for $r\geq 2$ even, so that $L\in  W^r_{\g-1}
 - W^{r+1}_{\g-1}$. From the hypothesis, we get that
\begin{equation}\label{hypo}
\dim T_L W^r_{\g-1} \cap T_LP^- \ge g-3={\rm dim}(P^-)-2.
\end{equation}
The Zariski tangent space to $W^r_{\g-1}$ at $L$ is given by the
orthogonal complement to the image of the Petri map (e.g.
\cite[Prop. 4.2]{ACGH}):
$$
H^0(\upC,L)\otimes H^0(\upC,\sigma^*(L))\rightarrow H^0(\upC,\omega_{\upC}),
$$
where we have used that $\omega_{\upC}=\pi^*(\omega_C)=L\otimes \sigma^*(L)$.
On the other hand, the tangent space to the Prym is by definition
$T_LP^-=H^0(\upC,\omega_{\upC})^-$, the $(-1)$-eigenspace of $H^0(\upC, \omega_{\upC})$
relative to the involution $\sigma$. Therefore, it is easy to see that
the intersection of the Zariski tangent spaces
$T_L W^r_{\g-1}\cap T_LP^-$ is given as the orthogonal
complement to the image of the map$$
v_0:\wedge^2H^0(\upC,L)\rightarrow H^0(\upC,\omega_{\upC})^-$$ defined
by $v_0(s_i\wedge s_j)=s_i\sigma^*s_j-s_j\sigma^*s_i.$

The inequality (\ref{hypo}) is equivalent to $\codim(\ker v_0)\le
2$. On the other hand, the decomposable forms in
$\wedge^2H^0(\upC,L)$ form a subvariety of dimension $2r-1\geq 3$, 
and so there is a decomposable vector $s_i\wedge s_j$ in $\ker
v_0$. This means that $s_i\sigma^*s_j-s_j\sigma^*s_i=0$, or in
other words that $\tfrac{s_j}{s_i}$ defines a rational function
$h$ in $C$. We conclude by taking $M=\O_C((h)_0)$ and $F$ the be
the maximal common divisor between $(s_i)_0$ and $(s_j)_0$.
\end{proof}

\begin{proof}[Proof of Theorem \ref{dim}]
The dimension of $ V(\upC)=V^2$ is at least $g-4$ by the theorem
of Bertram (\cite{B}, see also \cite{DP}). Suppose, by
contradiction, that there is an irreducible component $Z\subseteq
V^2$ such that $\dim Z=m\ge g-3$. Then, by applying the preceding
Lemma \ref{big-comp} for the general element $L\in Z$,
$$
L\cong \pi^*M\otimes \mathcal O_{\upC}(B)
$$
where
$M$ is an invertible sheaf on $C$ such that $h^0(M)\geq 2$, and $B$ is an effective divisor on
$\upC$ such that
$Nm(B)\in |K_C\otimes M^{\otimes -2}|$.
The family of such pairs $(M,B)$ is a finite cover of the set of pairs $\set{M, F}$
where:
\begin{itemize}
\item $M$ is an invertible sheaf on $C$ of degree $d\geq 2$ such that $h^0(M)\geq 2$,
\item $F$ is an effective divisor on $C$ of degree $2g-2-2d\geq 0$, such that
$F\in |K_C\otimes M^{\otimes -2}|$.
\end{itemize}

\noindent By Marten's theorem applied to the non-hyperelliptic curve $C$
(see \cite[Pag. 192]{ACGH}),
the dimension of the above family of line bundles $M$
is bounded above by
\begin{equation}\label{dis1}
{\rm dim}(W^1_d)< d-2.
\end{equation}

Fixing a line bundle $M$ as above, the dimension of possible $F$ satisfying the second condition is bounded by Clifford's theorem,
\begin{equation}\label{dis2}
h^0(K_C\otimes M^{\otimes -2})-1 \leq g-1-d,
\end{equation}

By putting together the inequalities (\ref{dis1}) and
(\ref{dis2}), we get that the dimension $m$ of our family of pairs
$\set{M, F}$ is bounded above by $m< d+2 +g-1-d=g-3$,
contradicting our hypothesis.
\end{proof}

\section{Computing the $V^i(\mathcal{I}_{\upC}(\Xi) )$}

In this section we compute the cohomological support loci for the ideal sheaf
$\I_{\upC}(\Xi)$, which can be identified
with the auxiliary canonical loci
$$\V^i(\I_{\upC}(\Xi^+))\subset P^-$$
defined in (\ref{auxiliary}).

\begin{teo}\label{1-twist}
The cohomological support loci for $\I_{\upC}(\Xi)$ are 
$$\begin{sis}
&V^0(\I_{\upC} (\Xi))=V^1(\I_{\upC} (\Xi))=V(\upC), \\
&V^2(\I_{\upC} (\Xi))=P^-,\\
&V^{\geq 3}(\I_{\upC} (\Xi))=\emptyset.
\end{sis}$$
\end{teo}
\begin{proof}
The equality $V^0(\I_{\upC} (\Xi))=V(\upC)$ is just the
definition of the theta-dual of $\upC$. Consider the exact
sequence defining the ideal sheaf $\I_{\upC}$ twisted by the
divisor $\Xi^+_{E}$, for $E\in P^-$:
\begin{equation*}
0\to \I_{\upC}(\Xi^+_{E})\to \O_{P'}(\Xi^+_{E})\to \O_{\upC}(\Xi^+_{E})\to 0.
\end{equation*}
By taking cohomology and using the vanishing
$H^{j}(P',\O_{P'}(\Xi^+_{E}))=0$ for $j>0$, we get the emptiness
of $\V^i(\mathcal{I}_{\upC}(\Xi^+) )$ for $i\geq 3$ and the two
exact sequences
$$0\ra H^0(\mathcal{I}_{\upC} (\Xi^+_{E}))\ra  H^0(\O_{P'}(\Xi^+_{E}))\stackrel{\psi_E}\ra H^0(\O_{\upC}(\Xi^+_{E}))\ra H^1(\I_{\upC}(\Xi^+_{E}))\ra 0,$$
$$0\ra  H^1(\O_{\upC}(\Xi^+_{E}))\ra H^2(\I_{\upC} (\Xi^+_{E}))\ra 0.$$
Using the second exact sequence and the Lemma \ref{restriction}, we get that
$$\V^2(\I_{\upC} (\Xi^+))=\set{E \mid h^1(\upC,E)>0}.
$$
Since $E$ has degree $\g-1$, by Riemann-Roch we have that
$h^1(\upC,E)= h^0(\upC,E)$, which is greater than $0$
for all $E\in P^-$ by the definition of $P^-$.

Consider now the first above exact sequence. Since
$H^0(P',\O_{P'}(\Xi^+_{E}))=1$ and $\V^1(\I_{\upC} (\Xi^+))$
consists of the elements $E$ such that the map $\psi_E$ is not
surjective, we get using again Lemmas \ref{restriction} and
\ref{dual-BN}
$$
\V^1(\I_{\upC} (\Xi^+))= \set{E \mid h^1(\upC,\O_{\upC} (\Xi^+_{E})
)>1 }\cup \set{E\mid h^0(P',\I_{\upC} (\Xi^+_{E}) )=1}=
$$
$$=\V^2 \cup V(\upC)=V(\upC).
$$
\end{proof}

\section{Computing the $V^i(\mathcal{I}_{C}(2\Xi) )$}

In this section we determine the cohomological support loci of the ideal sheaf $\I_{\upC}(2\Xi^+)$.  In the proof of the next Theorem, we will need to know that the set $S(\upC)$ of theta-characteristics of $\upC$ has a point in $P^-$.  It is not much more work to show the stronger lemma below.

 \begin{lema}\label{theta-char} A theta-characteristic $L\in S(\upC)$ belongs to ${\rm Nm}^{-1}(\omega_C)=P^+\cup P^-$ if and only if it is the pull-back of a theta-characteristic $M$ on $C$. Moreover, it holds that
 $$|S(\upC)\cap P^-|=|S(\upC)\cap P^+|=2^{2g-1}.$$ \end{lema} \begin{proof} If $M$ is a theta-characteristic on $C$, then $\pi^*(M)\in S(\upC)\cap {\rm Nm}^{-1}(\omega_C)$ as $$\begin{sis} &{\rm Nm}(\pi^*(M))=M^{\otimes 2}=\omega_C,\\ &\pi^*(M)^{\otimes 2}=\pi^*(M^{\otimes 2})=\pi^*(\omega_C)=\omega_{\upC}. \end{sis}$$ Conversely, if $L^{\otimes 2}=\omega_{\upC}$ and ${\rm Nm}(L)= \omega_C$, then $$ L\otimes L=\omega_{\upC}=\pi^*(\omega_C)=\pi^*({\rm Nm}(L))=L\otimes \sigma^*(L), $$ which implies that $L=\pi^*(M)$ for some $M\in \Pic(C)$. By applying the Norm map, we obtain $\omega_C={\rm Nm}(L)={\rm Nm}(\pi^*(M))=M^{\otimes 2}$.

Moreover, if we denote by $\eta_0$ the line bundle of order two on $C$ satisfying $\pi_*(\O_{\upC})=\O_C\oplus \eta_0$, then the pull-back $\pi^*(M)$ of a theta-characteristic $M$ on $C$ belongs to $P^-$ if and only if $h^0(C, M)\not\equiv h^0(C, M\otimes \eta_0) \mod 2$ (and to $P^+$ otherwise), as follows from the formula $H^0(\upC, \pi^*(M))=H^0(C, M)\oplus H^0(C, M\otimes \eta_0)$.

  Now, fix a theta-characteristic $M_0$ of $C$. Then all the theta-characteristics of $C$ are of the form $M_0\otimes \eta$ for a unique $\eta\in J_2(C)$, where $J_2(C)$ is the group of the $2^{2g}$ line bundles of $C$ whose square is trivial. Consider the following map $$\apl{q_0}{J_2(C)}{\Z/2 \Z}{\eta}{h^0(M_0\otimes \eta\otimes \eta_0)- h^0(M_0\otimes \eta)\mod 2} $$ The Riemann-Mumford relation (see \cite[p. 182]{M71}) yields $$q_0(\eta)=q_0(\eta_0)+{\rm ln}\,e_2(\eta, \eta_0),$$ where $e_2:J_2(C)\times J(C)\to \{\pm 1\}$ is the Riemann skew-symmetric bilinear form (recall that ${\rm char}(k)\neq 2$) and ${\rm ln}$ is defined by ${\rm ln}(+1)=0$ and ${\rm ln}(-1)=1$. From the nondegeneracy of the form $e_2$ and the fact that $\eta_0\neq \O_C$, it follows that the function $$\apli{J_2(C)}{\{\pm 1\}}{\eta}{e_2(\eta_0, \eta)}$$ takes half the times the value $+1$ and half the value $-1$. This concludes our proof. \end{proof}

To complete the main theorem announced in the introduction we have
the following:
\begin{teo}\label{2-twist}
For a fixed $E_0\in P^-$, the cohomological support loci for the ideal sheaf
$\I_{\upC}(2\Xi^+)$ are isomorphic to
$$ \begin{sis}
&V^0(\I_{\upC} (2\Xi^+))=
\begin{cases}
(P^-)^{\otimes 2} & \text{ if } g\geq 4,\\
\{\omega_{\upC}\}& \text{ if } g=3,
\end{cases}\\
&V^1(\I_{\upC} (2\Xi^+))=V^2(\I_{\upC} (2\Xi^+))=
\{\omega_{\upC}\},\\
&V^{\geq 3}(\I_{\upC} (2\Xi^+))=\emptyset.
\end{sis}$$
\end{teo}
\begin{proof}
We first compute the auxiliary canonical loci $\V^i(\I_{\upC}(2\Xi^+))\subset P^-$
introduced in (\ref{auxiliary}). To this end,
consider the exact sequence defining the ideal sheaf $\I_{\upC}$
twisted by the divisor $2\cdot \Xi^+_{E}$ with $E\in P^-$:
\begin{equation*}
0\to \I_{\upC}(2\cdot \Xi^+_{E})\to \O_{P'}(2\cdot \Xi^+_{E})\to
\O_{\upC}(2\cdot \Xi^+_{E})\to 0.
\end{equation*}
By taking cohomology and using the vanishing
$H^{j}(P',\O_{P'}(2\cdot \Xi^+_{E}))=0$ for $j>0$, we get the
emptiness of $\V^i(\mathcal{I}_{\upC}(2\Xi^+) )$ for $i\geq 3$ and
the two exact sequences
$$H^0(\mathcal{I}_{\upC} (2\cdot \Xi^+_{E}))\hookrightarrow  H^0(\O_{P'}(2\cdot \Xi^+_{E}))
\stackrel{\phi_{E}}\ra H^0(\O_{\upC}(2\cdot \Xi^+_{E}))\twoheadrightarrow
H^1(\I_{\upC}(2\cdot \Xi^+_{E})),$$
$$0\ra  H^1(\O_{\upC}(2\cdot \Xi^+_{E}))\ra H^2(\I_{\upC} (2\cdot \Xi^+_{E}))\ra 0.$$
By Lemma \ref{restriction}, we have that $\O_{\upC}(2\cdot \Xi^+_{E})=E^{\otimes 2}$
and therefore the second exact sequence implies that
$$\V^2(\I_{\upC} (2\Xi^+))=\set{E^{\otimes 2}\: | \: h^1(\upC,E^{\otimes 2})>0}=
S(\upC)\cap P^-.
$$
Moreover, from the Proposition \ref{surjectivity} below and the first exact sequence,
we have that
$$\V^1(\I_{\upC} (2\Xi^+))\subset S(\upC)\cap P^-.$$

In order to determine $\V^0(\I_{\upC} (2\Xi^+))$, we consider the first above exact sequence.
If $g\geq 4$, then we have the inequality
$$h^0(P',\O_{P'}(2\cdot \Xi^+_{E}))=2^{g-1}>2g-1=\g\geq
h^0(\upC,\O_{\upC}(2\cdot \Xi^+_{E})),$$
from which we conclude that
$$\V^0(\I_{\upC} (2\Xi^+))=P^- \text{ if } g\geq 4.$$
On the other hand, if $g=3$ and $E$ is not a theta characteristic, then
the map $\phi_{E}$ is a surjection (as we will prove in the Proposition \ref{surjectivity})
between two spaces of the same dimension and therefore an isomorphism.
This implies that
$$\V^0(\mathcal{I}_{\upC}(2\Xi^+))\subset
S(\upC)\cap P^- \text{ if }g=3.$$
\noindent Observe that $S(\upC)\cap P^-\neq\emptyset$ according to Lemma \ref{theta-char}.
Therefore, using our canonical identification (\ref{canonical-support}), we get for
the cohomological support loci
 $$ \begin{sis}
&V^0(\I_{\upC} (2\Xi^+))
\begin{cases}
= (P^-)^{\otimes 2} & \text{ if } g\geq 4,\\
\subset \{\omega_{\upC}\}& \text{ if } g=3,
\end{cases}\\
& V^1(\I_{\upC} (2\Xi^+))\subset \{\omega_{\upC}\},\\
&V^2(\I_{\upC} (2\Xi^+))=
\{\omega_{\upC}\},\\
&V^{\geq 3}(\I_{\upC} (2\Xi^+))=\emptyset.
\end{sis}$$

\noindent We deduce that the ideal $\I_{\upC} (2\Xi^+)$ is GV-sheaf in the sense of
\cite[Def. 3.1]{PP06} and hence, by \cite[Prop. 3.13]{PP06}, we get that
$$\{\omega_{\upC}\}=V^2(\I_{\upC} (2\Xi^+))\subset
V^1(\I_{\upC} (2\Xi^+))\subset V^0(\I_{\upC} (2\Xi^+)),$$
which gives the desired conclusion.
\end{proof}


The next Proposition is due to Beauville;  a proof is given in  \cite[Lemma 2.4]{IvS}
for the case of genus four curves.  We give a proof here for the sake of the completeness.

\begin{propo}\label{surjectivity}
 If $E\in P^-$ is not a theta-characteristic, then the restriction map
$$\phi_{E}:H^0(P',\O_{P'}(2\cdot \Xi^+_{E}))\ra H^0(\upC,\O_{\upC}(2\cdot \Xi^+_{E}))$$
is surjective.
\end{propo}
\begin{proof}
We start by proving that a general $\alpha\in \Xi^+_{E}\subset P'$ satisfies
\begin{enumerate}
 \item $(-1)^*(\alpha)=\alpha^{-1} \not\in \Xi^+_{E}$,
 \item $|E\otimes \alpha|$ is a pencil without base points.
\end{enumerate}

The first assertion follows from the fact that $\Xi_{E}^+$ is not
symmetric. Indeed, $\Xi_{E}^+$ is symmetric if and only if $E$ is
a theta-characteristic (see \cite[Thm. 11.2.4]{BL}), which we have
excluded by hypothesis.

\vspace{0,4cm}

The fact that the complete linear series $|E\otimes \alpha|$ is a
pencil for a general $\alpha \in \Xi^+_{E}$ follows from the fact
that a generic element $L\in \Xi^+$ has $h^0(L)=2$.

\vspace{0,4cm}

Now we want to see that the linear system $\abs{E\otimes \al}$ has no base points for a general $\al \in \Xi^+_E$. Consider the incidence variety
$$\xymatrix{&\Xi^+_E\times \upC\\
&I=\set{(\al,p)\mid p\text{ is a base point of }|\al\otimes E|}
\ar@{^{(}->}[u]\ar[ld]^{p_1}\ar[rd]_{p_2}\\
\Xi^+_E&&\upC}$$
For every point $q\in\upC$, the following injection
$$\apli{p_2^{-1}(q)}{V^2}{\alpha}{\alpha\otimes E\otimes\O(-q+\sigma q)}$$
is well-defined since $q$ is a fixed point of $|E\otimes \alpha|$.
Therefore, by Theorem \ref{dim}, the fibers of $p_2$ have
dimension at most $g-4$ and hence $I$ has dimension at most $g-3$.
Since the dimension of $\Xi^+_{E}$ is $g-2$, this implies that the
first projection is not dominant and hence the conclusion.

\vspace{0,4cm}

Now we want to find elements in $H^0(P',\O_{P'}(2\cdot \Xi^+_{E}))$ that form a basis when restricted to $H^0(\upC,\O_{\upC}(2\cdot \Xi^+_{E}))$. From Lemma \ref{restriction} and the fact that
$E$ is not a theta characteristic, we get
that $h^0(\upC,\O_{\upC}(2\cdot \Xi^+_{E})=h^0(\upC, E^{\otimes 2})=\g-1$.
We begin by choosing an $\alpha\in \Xi^+_{E}\subset P^-$ satisfying the two conditions above.
In particular, from the condition $(2)$, we can choose an effective divisor
$E_\al=\sum_{1\leq i\leq \g-1} x_i\in |E\otimes \alpha|$ such that all the points $x_i$
are distinct.

\vspace{0,4cm}


We define the following effective divisors in $P^-$
$$E_{\al,j}:= E_\al-x_j+\sigma x_j=\sigma x_j +\sum_{i\neq j}x_i \qquad \text{ for }j=1,\ldots,\g-1.$$
By condition $(2)$, we get that $h^0(P',\O_{P'}(E_{\al,j}))=1$ and therefore, by
Lemma \ref{restriction}, $\left(\Xi^+_{E_{\al,j}}\right)_{\mid \upC}= E_{\al,j}$.

\vspace{0,4cm}

Consider next the line bundle $E\otimes \al^{-1}\otimes
\O_{\upC}(x_j-\sigma x_j)\in P^-$. Since $h^0(\upC,E\otimes
\al^{-1})=0$ by condition $(1)$, using Mumford's parity trick we
deduce that
$$h^0(\upC,E\otimes \al^{-1}\otimes \O_{\upC}(x_j-\sigma x_j))=1.$$
Define $E'_{\al,j}$ to be the unique effective divisor of
$|E\otimes \al^{-1}\otimes \O_{\upC}(x_j-\sigma x_j)|$
By Lemma \ref{restriction}, we get that $\left(\Xi^+_{E'_{\al,j}}\right)_{\mid\upC}= E'_{\al,j}$.

\vspace{0,4cm}

Summing up, we have constructed $\g-1$ couples of divisors
$(E_{\alpha,j},E'_{\alpha,j})$ satisfying
$$\begin{sis}
&\O_{P'}(\Xi^+_{E_{\al,j}}+\Xi^+_{E'_{\al,j}})\cong \O_{P'}(2\cdot \Xi^+_{E}),\\
&(\Xi^+_{E_{\al,j}}+\Xi^+_{E'_{\al,j}})_{\mid\upC}=E_{\al,j}+E'_{\al,j}.
\end{sis}$$
It remains to show that the $\g-1$ divisors $ E_{\al,j}+E'_{\al,j} $ corresponds to independent
sections on $H^0(\upC, 2\cdot E)$. This will follow from the next

\vspace{0,4cm}

\noindent{\bf Claim}: $x_j\not \in E_{\alpha,j}+E'_{\alpha, j}$
and $x_j\in E_{\alpha,k}+E'_{\alpha, k}$ for every $k\neq j$.\\
\noindent By the definition of the $E_{\alpha,j}$ and using that
$\sigma$ has no fixed points, we get that $x_j\not\in
E_{\alpha,j}$ while $x_j\in E_{\alpha,k}$ for every $k\neq j$.
Finally, observe that $\O_{\upC}(E'_{\alpha,j})\otimes
\O_{\upC}(\sigma x_j-x_j)= E\otimes \alpha^{-1}$ which by
condition $(1)$ has no sections. This can happen only if
$E'_{\alpha,j}-x_j$ is not effective, or in other words $x_j\not
\in E'_{\alpha,j}$.

\end{proof}

\end{document}